\documentclass{amsart}

\newtheorem{theorem}{Theorem}[section]

\newtheorem{corollary}[theorem]{Corollary}

\begin{document}

\title[Isotropic hyperelastic constitutive laws]{On the dual variable of the Cauchy stress tensor\\
in isotropic finite hyperelasticity}



\author{Claude VALLEE}


\address{Laboratoire de M\'{e}canique des Solides, UMR CNRS 6610, Universit\'{e} de Poitiers, SP2MI, T\'{e}l\'{e}port 2,
Boulevard Marie et Pierre Curie, B.P. 30179, 86962, Futuroscope-Chasseneuil Cedex, France \\
              Tel.: 0033-(0)549-496798\\
              Fax: 0033-(0)549-496791}
\email{vallee@lms.univ-poitiers.fr}           

\author{Danielle FORTUNE}
\address{Laboratoire de M\'{e}canique des Solides, UMR CNRS 6610, Universit\'{e} de Poitiers, SP2MI, T\'{e}l\'{e}port 2,
Boulevard Marie et Pierre Curie, B.P. 30179, 86962, Futuroscope-Chasseneuil 
Cedex, France} 
\author{Camelia LERINTIU} 
\address{Laboratoire de M\'{e}canique des Solides, UMR CNRS 6610, Universit\'{e} de Poitiers, SP2MI, T\'{e}l\'{e}port 2,
Boulevard Marie et Pierre Curie, B.P. 30179, 86962, Futuroscope-Chasseneuil 
Cedex, France}



\begin{abstract}
Elastic materials are governed by a constitutive law linking the second Piola-Kirchhoff stress tensor $\Sigma$ and the right Cauchy-Green strain tensor $C=F^TF$. Isotropic elastic materials are the special ones for which the Cauchy stress tensor $\sigma$ depends solely of the left Cauchy-Green strain tensor $B=FF^T$. In this paper we revisit the following property of isotropic hyperelastic materials: if the constitutive law linking $\Sigma$ and $C$ derives from a potential $\alpha$, then $\sigma$ and $\ln B$ are linked by a constitutive law deriving from the potential $\alpha \circ \exp$. We give a new and concise proof which is based on an explicit formula expressing the derivative of the exponential of a tensor.
\end{abstract}
\maketitle

\section{Introduction}
\label{intro}
According to the mass conservation principle, the mass density per unit volume $\rho$ and its initial value $\rho_0$ are in the ratio
$$\frac{\rho_0}{\rho}=\det F=(\det C)^\frac{1}{2}=(\det B)^\frac{1}{2}.$$
The relation 
$$(\det F)\sigma=F\Sigma F^T$$
between the Cauchy stress tensor $\sigma$ and the second Piola-Kirchhoff stress tensor $\Sigma$ can be rewritten
$$\frac{\sigma}{\rho}=F\frac{\Sigma}{\rho_0}F^T.$$
Let us agree to formulate the elastic materials constitutive laws as: 
$$\frac{\Sigma}{\rho_0}=h(C).$$
The polar decomposition $F=RU$ of the deformation gradient \cite{Cia88} implies:
$$B=FF^T=RU^2R^T=RCR^T \textnormal{ or  } C=R^TBR.$$ 
This allows to translate the relation between $\Sigma$ and $C$ by a law satisfied by $\sigma$: 
$$\frac{\sigma}{\rho}=RUh(R^TBR)UR^T=(RUR^T)Rh(R^TBR)R^T(RUR^T)$$
where we have enlightened the tensor $RUR^T$ which is nothing else than the square root $B^\frac{1}{2}$ of the positive definite symmetric tensor $B$. A priori, for elastics materials, the tensor $\frac{\sigma}{\rho}$ is a function of $B$ and $R$: 
$$\frac{\sigma}{\rho}=B^\frac{1}{2}Rh(R^TBR)R^TB^\frac{1}{2}.$$
It will depend solely of $B$ in a single case: when the tensor $Rh(R^TBR)R^T$ does not depend on the rotation $R$. The rotations forming a group, the only possible tensorial functions $h$ are those satisfying the relations of isotropy with respect to $B$:
$$Rh(R^TBR)R^T=h(B)\quad \textnormal{or}\quad R^Th(B)R=h(R^TBR).$$
Because of the relation $C=R^TBR$, the isotropy of  the function $h$ can alternatively be expressed with respect to $C$: 
$$Rh(C)R^T=h(RCR^T)\quad \textnormal{or}\quad R^Th(RCR^T)R=h(C).$$
To summarize: \textit{if the law $\frac{\Sigma}{\rho_0}=h(C)$ is isotropic, then $\frac{\sigma}{\rho}$ depends only of $B$, and it is the sole case ; furthermore, under this isotropy condition} 
$$\frac{\sigma}{\rho}=B^\frac{1}{2}h(B)B^\frac{1}{2}.$$
In this paper, we revisit the property of isotropic hyperelastic materials for which the existence of a potential expressing the constitutive law between $\frac{\Sigma}{\rho_0}$ and $C$ implies the existence of a potential linking $\frac{\sigma}{\rho}$ and $\ln B$. 
\section{Isotropy of the constitutive law linking $\frac{\sigma}{\rho}$ and $B$}
\label{sec:1}
Let $\Omega$ be a rotation, if we change $B$ into $\Omega^T B \Omega$, then $B^\frac{1}{2}$ is changed into $\Omega^T B^\frac{1}{2} \Omega$ and $\frac{\sigma}{\rho}$ is changed in: 
$$\Omega^T B^\frac{1}{2} \Omega h(\Omega^TB\Omega) \Omega^T B^\frac{1}{2} \Omega=\Omega^T B^\frac{1}{2} h(B) B^\frac{1}{2} \Omega=\Omega^T  \frac{\sigma}{\rho} \Omega.$$
The isotropy of the constitutive law linking $\frac{\Sigma}{\rho_0}$ and $C$ is thus transferred to the constitutive law linking $\frac{\sigma}{\rho}$ and $B$. 
\section{Coaxiality of $B$ and $h(B)$}
\label{sec:2}
%
%
\begin{theorem} Because $h$ is isotropic, the symmetric tensors $B$ and $h(B)$ are coaxial (i.e. they have the same eigenvectors).
\end{theorem}
\begin{proof}
Let $n$ be an eigenvector of $B$ chosen unitary, and let us consider the rotation of angle $\pi$ around $n$: 
$$S=(\cos \pi)I+(1-\cos \pi)nn^T=2nn^T-I$$
with $I$ as the identity tensor. Such a symmetry $S$ leaves $n$ unchanged and changes any orthogonal vector to $n$ in its opposite. The tensor $B$ being symmetric, its other two eigenvectors are orthogonal to $n$, as a consequence $S^TBS=B$. \\
The isotropy condition implies $S^Th(B)S=h(S^TBS)$ or $h(B)S=Sh(B)$, therefore $h(B)[Sn]=S[h(B)n]$ or $S[h(B)n]=h(B)n$. Since the sole vectors unchanged by $S$ are the vectors parallel to $n$, the last equality is possible only when the vector $h(B)n$ remains parallel to the vector $n$, that is to say when $n$ is also an eigenvector for $h(B)$. 
\end{proof}
We easily deduce from this coaxiality property the two next corollaries, which will reveal important in the following.
%
%
\begin{corollary} $B$ and $h(B)$ commute.
\end{corollary}
%
%
\begin{corollary} For every real number $s$, $h(B)$ commutes with the power $B^s$ of $B$. 
\end{corollary} 
The choice $s=\frac{1}{2}$ allows one to simplify the expression $\frac{\sigma}{\rho}=B^\frac{1}{2} h(B) B^\frac{1}{2}$ in
$$\frac{\sigma}{\rho}=h(B)B.$$ 
\section{Hyperelastic materials}
\label{sec:3}
\subsection{Existence of a potential between the second Piola-Kirchhoff stress tensor $\Sigma$ and the right Cauchy-Green strain tensor $C$}
\label{sec:4}
Let us consider a derivable function $\alpha$ of $C$, its derivative $D\alpha(C)$ is a linear mapping from the space of symmetric tensors to $\mathbb{R}$. Thus, there exists a symmetric tensor denoted $\frac{\partial \alpha}{\partial C}$ such that for every variation $\delta C$ of $C$: 
$$D\alpha(C)\delta C=\mathrm{tr}(\frac{\partial \alpha}{\partial C}\delta C).$$
Hyperelastic materials are those for which there exists a function $\alpha$ such that
$$\frac{\Sigma}{\rho_0}=\frac{\partial \alpha}{\partial C}$$
In this assumption, we will say that the constitutive law linking the tensors $\frac{\Sigma}{\rho_0}$ and $C$ derives from the potential $\alpha$. 
\subsection{Derivative of the exponential of a matrix}
\label{sec:5}
Let us consider a square matrix $A$ and a real number $t$, the exponential $\exp(tA)$ is the solution of the matricial ordinary differential equation 
$$\frac{d}{dt}\exp(tA)=A\exp(tA)$$
which is equal to $I$ at $t=0$. Let $\delta A$ be a variation of $A$, in the varied equation 
$$\frac{d}{dt}D(\exp)(tA)(t\delta A)=\delta A[\exp(tA)]+AD(\exp)(tA)(t\delta A)$$
let us introduce the square matrix $M(t)$ defined by 
$$D(\exp)(tA)(t\delta A)=[\exp(tA)]M(t)$$
The varied equation becomes
$$[\frac{d}{dt}\exp(tA)]M(t)+[\exp(tA)]\frac{dM}{dt}=\delta A[\exp(tA)]+A[\exp(tA)]M(t)$$
and simplifies itself into the ordinary differential equation
$$\frac{dM}{dt}=[\exp(-tA)]\delta A[\exp(tA)]$$
which can be integrated by quadrature. Because $M(0)$ vanishes, we easily deduce from it the value of $M(1)$ and thereafter the variation of the exponential of a matrix \cite{Sou59}: 
$$D(\exp)(A)(\delta A)=[\exp(A)]\int_0^{1}[\exp(-sA)]\delta A [\exp(sA)]ds.$$
In the special case where $A$ is the logarithm of the positive definite tensor $B$, this formula allows us to predict that for every variation $\delta B$ of $B$: 
$$D(\exp)(\ln B)\delta B=B\int_0^{1}B^{-s}\delta B B^{s}ds$$
%
%
\subsection{Existence of a potential between the Cauchy stress tensor and the logarithm of the left Cauchy-Green strain tensor}
\label{sec:6}
%
%
%
\begin{theorem} If the tensor $\frac{\Sigma}{\rho_0}$ derives from a potential $\alpha$ of the tensor $C$, then the tensor $\frac{\sigma}{\rho}$ derives from the potential $\alpha \circ \exp$ of the tensor $\ln B$.
\end{theorem}
\begin{proof}
By deriving the compound function $\alpha \circ \exp$, we find successively: 
$$D(\alpha \circ \exp)(\ln B)\delta B=D\alpha(B)(D(\exp)(\ln B)\delta B)$$
$$=\mathrm{tr}(\frac{\partial \alpha}{\partial B}[D(\exp)(\ln B)\delta B])=\mathrm{tr}(h(B)B\int_0^{1}B^{-s}\delta B B^{s}ds)=\int_0^{1}\mathrm{tr}[h(B)BB^{-s}\delta B B^{s}]ds.$$
To simplify the last integral, it is necessary to pay attention on the switchings because the matrix $\delta B$ does not commute with the others. However, under the trace, we can make cross at the beginning the last term of the product of 5 matrices. Then from Corollary 2, we can switch this term $B^s$ with $h(B)$ and afterwards with $B$, it ends up just before $B^{-s}$. The product of the two matrices $B^s$ and $B^{-s}$ reduces to the identity tensor $I$, and the integral simplifies itself into
$$\mathrm{tr}(h(B)B\delta B)=\mathrm{tr}(\frac{\sigma}{\rho}\delta B).$$
The final value of the integral allows to conclude to the constitutive law:
$$\frac{\sigma}{\rho}=\frac{\partial (\alpha \circ \exp)}{\partial (\ln B)}$$
\end{proof}
%
%
\section{Conclusion}
\label{sec:7}
Without resorting to the Taylor expansion of the logarithm \cite{Val78} or of the exponential \cite{San01} of a symmetric tensor, nor to its spectral decomposition \cite{Hog87}, we have given an intrinsic proof of the existence of the potential $\alpha \circ \exp$ between $\frac{\sigma}{\rho}$ and $\ln B$. Numerous isotropic hyperelastic constitutive laws expressing directly $\sigma$ in term of $\ln B$ have been proposed (\cite{Bru01}, \cite{Ogd04}, \cite{Per02}, \cite{Sen05}, \cite{Xia02}) and numerically implemented \cite{Fen06}.\\ 
When the potential $\alpha \circ \exp$ is convex, the consideration of its Legendre-Fenchel-Moreau transform is a tool to perform the inversion of the constitutive law (\cite{Blu92}, \cite{Xia03}, \cite{Xia04}), ie to express the Hencky logarithmic strain tensor $\ln B$ in term of the Cauchy stress tensor $\sigma$.
%
%
%

%

\begin{thebibliography}{ab}
%
%
\bibitem{Blu92} {Blume, J.A.}: On the form of the inverted stress-strain law for isotropic hyperelastic solids.
International Journal of Non-Linear Mechanics. {\bf 27(3)}, 413--421 (1992) 
%
%
\bibitem{Bru01} {Bruhns, O.T., Xiao, H., Meyers, A.}: Constitutive inequalities for an isotropic elastic strain-energy function based on Hencky's logarithmic strain tensor. Proceedings of the Royal Society of London, Series A-Mathematical Physical and Engineering Sciences. {\bf 457(2013)}, 2207--2226 (2001) 
%
%
\bibitem{Cia88} {Ciarlet, P.G.}: Mathematical Elasticity, vol.1. 3D Elasticity, North-Holland, Amsterdam (1988) 
%
%
\bibitem{Fen06} {Feng, Z.-Q., Vall\'{e}e, C., Fortun\'{e}, D., Peyraut, F.}: The 3e hyperelastic model applied to the modeling of 3D impact problems. Finite Elements in Analysis and Design. {\bf 43(1)}, 51--58 (2006)
%
%
\bibitem{Hog87} {Hoger, A.}: The stress conjugate to logarithmic strain. International Journal of Solids and Structures. {\bf 23(12)}, 1645--1656 (1987) 
%
%
\bibitem{Ogd04} {Ogden, R.W., Saccomandi, G., Sgura, I.}: Fitting hyperelastic models to experimental data. Computational Mechanics. {\bf 34(6)}, 484--502 (2004)
%
%
\bibitem{Per02} {Peric, D., Owen, D.R.J., Honnor, M.E.}: A model for Finite Strain elastoplasticity based on logarithmic strains - computational issues. Computer Methods in Applied Mechanics and Engineering. {\bf 94(1)}, 35--61 (1992)
%
%
\bibitem{San01} 
{Sansour C.}: On the dual variable of the logarithmic strain tensor, the dual variable of the Cauchy stress tensor, and related issues. International Journal of Solids and Structures. {\bf 38(50-51)}, 9221--9232 (2001) 
%
%
\bibitem{Sen05} {Sendova, T., Walton, J.R.}: On strong ellipticity for isotropic hyperelastic materials based upon logarithmic strain. International Journal of Non-Linear Mechanics. {\bf 40(2-3)}, 195--212 (2005) 
%
%
\bibitem{Sou59}
{Souriau, J.M.}: Calcul lin\'eaire. P.U.F, Paris (1959)
%
%
\bibitem{Val78} 
{Vallé\'{e}e, C.}: Laws of isotropic hyperelastic behaviour (in french). International Journal of Engineering Science. {\bf 16(7)}, 451--457 (1978) 
%
%
\bibitem{Xia02} {Xiao, H., Chen, L.S.}: Hencky's elasticity model and linear stress-strain relations in isotropic finite hyperelasticity. Acta Mechanica. {\bf 157(1-4)}, 51--60 (2002)
%
%
\bibitem{Xia03} {Xiao, H., Chen, L.S.}: Hencky's logarithmic strain and dual stress-strain and strain-stress relations in isotropic finite hyperelasticity. International Journal of Solids and Structures. {\bf 40(6)}, 1455--1463 (2003)
%
%
\bibitem{Xia04} 
{Xiao, H., Bruhns, O.T., Meyers, A}: Explicit dual stress-strain and strain-stress relations of incompressible isotropic hyperelastic solids via deviatoric Hencky strain and Cauchy stress. Acta Mechanica. {\bf 168(1-2)}, 21--33 (2004)
%
\end{thebibliography}
\end{document}